\documentclass{owrart}

\newtheorem{Theorem}{Theorem}
\newtheorem{Proposition}[Theorem]{Proposition}
\newcommand{\op}[1]{\ensuremath{\mathrm{#1}}}
\newcommand{\ol}[1]{\ensuremath{\overline{#1}}}
\newcommand{\reg}{\op{reg}}

\newcommand{\fg}{\ensuremath{\mathfrak{g}}}
\newcommand{\dd}{\op{\tt{d}}}
\newcommand{\Ad}{\ensuremath{\operatorname{Ad}}}
\newcommand{\ad}{\ensuremath{\operatorname{ad}}}

\newcommand{\Aut}{\ensuremath{\operatorname{Aut}}}
\newcommand{\der}{\ensuremath{\operatorname{der}}}

\newcommand{\from}{\ensuremath{\nobreak\colon\nobreak}}
\renewcommand{\to}{\ensuremath{\nobreak\rightarrow\nobreak}}

\begin{document}


\begin{talk}[Christoph Wockel]{Bas Janssens}
{OWR: Universal central extensions of gauge groups}
{Janssens, Bas}

\noindent 
We indicate how to calculate the universal central extension of the gauge algebra
$\Gamma(\mathrm{ad}(P))$, and how to obtain from this the corresponding universal 
central extension of the gauge group $\Gamma(\mathrm{Ad}(P))$. 
 
Gauge groups occur as vertical symmetries of  
gauge theories, in which fields are connections on a principal 
$G$-bundle $P \rightarrow M$, and 
the action is 
invariant under vertical
automorphisms of $P$. 
If we set $\mathrm{Ad}(P) := P \times_{\mathrm{Ad}} G$
and similarly 
$\mathrm{ad}(P) := P \times_{\mathrm{ad}} \mathfrak{g}$ (with $\mathfrak{g}$
the Lie algebra of $G$), we can identify the group
of vertical automorphisms with $\Gamma(\mathrm{Ad}(P))$, and
its Lie algebra with $\Gamma(\mathrm{ad}(P))$.

In the case that $P$ admits a \emph{flat} equivariant connection, these gauge algebras
closely resemble
equivariant map algebras and (twisted multi) loop algebras.
Using the flat connection, one 
finds a cover $N \to M$, a monodromy group $\Delta < \pi_1(M)$, and
a homomorphism $\Delta \to G$ such that
$P = N \times_{\Delta} G$.
The adjoint bundle then takes the shape 
$\mathrm{ad}(P) = N \times_{\Delta}\mathfrak{g}$,
so that the gauge algebra is just  
$\Gamma(\mathrm{ad}(P)) = (C^{\infty}(N,
\mathbb{R})\otimes_{\mathbb{R}}\mathfrak{g})^{\Delta}$, 
the Lie algebra of smooth equivariant maps from $N$ to $\mathfrak{g}$.

If $X$ is an affine variety over $\mathbb{R}$ with an action of a discrete group 
$\Delta$, and 
$\Delta$ acts by automorphisms on a real Lie algebra $\fg$, then the equivariant 
map algebra $(\mathbb{C}[X] \otimes_{\mathbb{R}} \mathfrak{g})^{\Delta}$ is the 
Lie algebra of equivariant regular maps from $X$ to $\fg_{\mathbb{C}}$. 

The set
$X^{\mathrm{reg}}_{\mathbb{R}}$
of regular real points constitutes a smooth manifold,
and under suitable conditions (see prop.~\ref{density}), the homomorphism
$\mathbb{C}[X] \to C^{\infty}(X^{\mathrm{reg}}_{\mathbb{R}},\mathbb{C})$
is injective with dense image.
If the action of $\Delta$ restricts to
$X_{\mathbb{R}}^{\mathrm{reg}}$, then we obtain an inclusion
$(\mathbb{C}[X] \otimes_{\mathbb{R}} \mathfrak{g})^{\Delta} \hookrightarrow
(C^{\infty}(X_{\mathbb{R}}^{\mathrm{reg}}, 
\mathbb{C}) \otimes_{\mathbb{R}} \mathfrak{g})^{\Delta}$
of Lie algebras with dense image.
If moreover 
$X_{\mathbb{R}}^{\mathrm{reg}}/\Delta$
is a manifold, then we have realised the equivariant map algebra as a dense 
subalgebra of the complexification of the gauge algebra $\Gamma(\ad(P))$,
with $P \to X_{\mathbb{R}}^{\mathrm{reg}}/\Delta$ the principal 
$\mathrm{Aut}(\mathfrak{g})$-bundle 
$P = X_{\mathbb{R}}^{\mathrm{reg}} \times_{\Delta} \mathrm{Aut}(\mathfrak{g})$.

For example, let $X$ be
$T^{n} = \{(\vec{z},\vec{w}) \in \mathbb{C}^{2n} \,|\, 
z_k^{2} + w_k^{2} = 1, 1\leq k \leq n\}$, the complex
$n$-torus. 
In this case, 
$\mathbb{C}[T^{n}] \hookrightarrow C^{\infty}(T_{\mathbb{R}}^{n},\mathbb{C})$
is injective with dense image by Fourier theory.   
We now look for a regular group action on $T^{n}$ that restricts to 
$T_{\mathbb{R}}^{n}$ and such that $M = T_{\mathbb{R}}^{n} / \Delta$ is a manifold.
Although the Bieberbach groups spring to mind,
the choice that is studied most is 
$\Delta = \prod_{k=1}^{n} \mathbb{Z}/r_k\mathbb{Z}$,
with $\delta : (z_k \pm i w_k) \mapsto e^{\pm 2\pi i \delta_k / r_k} (z_k \pm i w_k)$.
In this case,
$M = T_{\mathbb{R}}^{n}/\Delta$ is again a torus. 
For any homomorphism $\Delta \to \mathrm{Aut}(\mathfrak{g})$, the
twisted multiloop algebra $(\mathbb{C}[T^{n}] \otimes_{\mathbb{R}}
 \mathfrak{g})^{\Delta}$
forms a dense subalgebra of the complexification of 
$\Gamma(\mathrm{ad}(P))$, where $P$ is the principal 
$\mathrm{Aut}(\mathfrak{g})$-bundle
$P = T_{\mathbb{R}}^{n} \times_{\Delta}\mathrm{Aut}(\mathfrak{g})$ over
$T^{n}_{\mathbb{R}}$.

The case of the circle is special in that every principal $G$-bundle over 
$M = T^{1}_{\mathbb{R}}$
is given by a twist $g\in G$ 
upon a full rotation, and therefore admits a flat connection. 
A smooth path connecting $g$ to $g'$ yields an isomorphism
of the corresponding bundles, so principal $G$-bundles are classified by $\pi_0 (G)$.
Consequently, the complexified adjoint bundles are classified by 
$\pi_0(\Aut(\fg_{\mathbb{C}}))$, which for simple $\fg_{\mathbb{C}}$
amounts to
diagram
automorphisms of order 1, 2 or 3.  
Complexified gauge algebras over $T^{1}_{\mathbb{R}}$ 
are thus precisely the closures of twisted loop algebras.


We return to the case of smooth principal fibre bundles which do not necessarily 
have a flat connection, and sketch the universal 2-cocycle
for the compactly 
supported gauge algebra $\Gamma_{\mathrm{c}}(\ad(P))$. 
We refer the interested reader to~\cite{bj_JW10} for details. 

For any Lie algebra $\fg$, 
set $V(\fg) := (\fg \otimes_s \fg) / \der(\fg) \cdot (\fg \otimes_s \fg)$,
and denote by $\kappa \colon \fg \times \fg \to V(\fg)\,;\, (x,y) 
\mapsto [x \otimes_{s}y]$ the universal $\der(\fg)$-invariant
bilinear form on~$\fg$. 
Any Lie connection $\nabla$ on $\ad(P)$ induces a flat connection
$\dd$ on the vector bundle $V(\ad(P)) \to M$, which does
not depend on $\nabla$ as any two Lie connections differ by
a pointwise derivation, which acts trivially on $V(\ad(P))$.
 Using the identities
 $\dd \kappa(\xi,\eta) = \kappa(\nabla\xi,\eta) + \kappa(\xi,\nabla\eta)$
 and 
 $\nabla [\xi,\eta] = [\nabla \xi , \eta] + [\xi, \nabla \eta]$
for all sections $\xi, \eta \in \Gamma_{\mathrm{c}}(\ad(P))$,
one checks that
 \begin{equation}\label{eqn:canonical-cocycle}
  \omega_{\nabla}\from \wedge^{2} \Gamma_{\mathrm{c}}(\ad(P))\to
  \ol{\Omega}{}^{1}_{\mathrm{c}}(M,V(\ad(P))
  \quad \xi \wedge \eta\mapsto 
  [\kappa(\xi,\nabla \eta)]
 \end{equation}
 defines a Lie algebra cocycle, where 
 the subscript $\mathrm{c}$ denotes compact support, and we set
$\ol{\Omega}{}^{1}_{\mathrm{c}} :=
 \Omega{}^{1}_{\mathrm{c}} / \dd \Omega^{0}_{\mathrm{c}}$.
 If $\fg$ is semisimple, then the cohomology class $[\omega_{\nabla}]$ does
 not depend on $\nabla$. We equip
 our
 spaces of smooth forms and sections with the usual LF-topology, 
in terms of which the 
universality result is 
formulated as follows.
 \begin{Proposition}\label{prop:universal_algebra_cocycle}
If $\fg$ is semisimple, then $[\omega_{\nabla}]$ is universal;
every continuous 2-cocycle $\psi$ with values in a trivial 
real topological module $W$ can be written up to coboundary
as $\psi = \phi \circ \omega_{\nabla}$, for some 
continuous $\mathbb{R}$-linear
$\phi : \ol{\Omega}{}^{1}_{\mathrm{c}}(M,V(\ad(P))) \to W$.
%
 \end{Proposition}
In \cite{bj_JW10},
this is proved by noting that 2-cocycles are automatically diagonal, 
so that the second cohomology in fact constitutes a sheaf.
The result can then be reduced to the well known local one, 
described e.g.\ in \cite{bj_Maier02}. 
Using the results of \cite{bj_Neeb02a,bj_Neeb02b,bj_NeebWockel}, 
this can be used (see \cite{bj_JW10}) to prove the following theorem.
\begin{Theorem}\label{thm:universal_group_extension}
 Let $P \to M$ be a principal fibre bundle with compact connected base,
 and with a semisimple structure group with finitely many connected components.
 Then the cocycle (\ref{eqn:canonical-cocycle}) integrates to a central
 extension of $\Gamma(\Ad(P))$ that is universal for abelian Lie groups 
modelled on Mackey-complete
locally convex spaces.
 \end{Theorem} 
 
Although the application of differential geometric techniques in an
algebraic context has intrinsic drawbacks, 
%
it is perhaps worth while to briefly explore the ramifications of  
proposition \ref{prop:universal_algebra_cocycle} to equivariant map algebras. 
We start by substantiating our claim 
as to the injectivity and denseness of  
$\mathbb{C}[X] \rightarrow C^{\infty}(X_{\mathbb{R}}^{\mathrm{reg}},\mathbb{C})$. 

\begin{Proposition}\label{density}
Let $X$ be an affine variety over $\mathbb{R}$ such that every connected component of
$X^{\mathrm{an}}$ possesses a nonsingular real point.
Then the ring homomorphism
$\mathbb{C}[X] \hookrightarrow C^{\infty}(X_{\mathbb{R}}^{\mathrm{reg}},\mathbb{C})$
is injective, with
dense image in the topology of uniform convergence of derivatives
on compact subsets.
\end{Proposition}
\begin{proof}
First, we prove that the image is dense.
As every smooth function on $X_{\mathbb{R}}^{\mathrm{reg}}$ can be approximated by
compactly supported smooth functions, and every compactly supported 
(in $X_{\mathbb{R}}^{\mathrm{reg}}$)
smooth function on $X_{\mathbb{R}}^{\mathrm{reg}}$ extends to a 
compactly supported (in $\mathbb{R}^{n}$)
smooth function on $\mathbb{R}^{n}$, it is enough to show
that every smooth function $f$ on a compact subset $K$ of $\mathbb{R}^{n}$
can be approximated by polynomials.
Now by Weierstra\ss' theorem, there exist,
for any multi-index $\vec{\mu}$, 
polynomials $p$ with $\sup_{K} |\partial_{\vec{\mu}} f - p|$ arbitrarily
small.  
By integrating these, one can produce polynomials $p$ such that 
$\sup_{K} |\partial_{\vec{\nu}} f - \partial_{\vec{\nu}}p|$
is arbitrarily small for all $\vec{\nu} < \vec{\mu}$.
A sequence $p_k$ of such polynomials for  $\vec{\mu}_{k} \to \infty$
(in the sense that for every fixed $\vec{\nu}$, we eventually have $\vec{\nu} < \vec{\mu}_{k}$)
will converge to $f$ uniformly on $K$ for every derivative. 

Next, we prove injectivity.
Denote by $\mathcal{O}^{\mathrm{an}}_Y$ and $C^{\infty}_Y$
the sheaves of analytic and smooth functions on $Y$.
Choose a nonsingular real point $x_i$ in each connected component (in the analytic topology)
of $X^{\mathrm{an}}$, so that 
$\mathbb{C}[X] \rightarrow \bigoplus_i 
\mathcal{O}_{X, x_i}^{\mathrm{an}}$
is injective.
Using the inverse function
theorem, we find analytic
charts $\phi_i : \mathbb{C}^{d} \supset U_i \to V_i \subset X^{\mathrm{an}}$ 
around $x_i$ 
in which
$U_i \cap \mathbb{R}^{d}$ corresponds to $V_i \cap X^{\reg}_{\mathbb{R}}$.
In those coordinates,  the map
$\mathcal{O}_{X, x_i}^\mathrm{an} \rightarrow 
C^{\infty}_{X^{\reg}_\mathbb{R}, x_i}$ corresponds to the injective map
$\mathcal{O}_{\mathbb{C}^{d} , 0}^\mathrm{an} \rightarrow 
C^{\infty}_{\mathbb{R}^{d} , 0}$, 
and is therefore injective.
Since the injective map $\mathbb{C}[X] \rightarrow \bigoplus_i 
C^{\infty}_{X^{\reg}_\mathbb{R}, x_i}$
factors through 
$\mathbb{C}[X] \rightarrow C^{\infty}(X^{\reg}_\mathbb{R},\mathbb{C})$,
the latter must be injective itself.\end{proof}

Consider $X$,  $\fg$ and $\Delta$ as before, but now with 
$X^{\mathrm{reg}}_{\mathbb{R}}/\Delta$ a \emph{compact} manifold,
and $\fg$ \emph{semisimple}. 
Since $ \iota \colon (\mathbb{C}[X]\otimes_{\mathbb{R}}\fg)^{\Delta} \hookrightarrow 
\Gamma(\ad(P)_{\mathbb{C}})$
is a dense inclusion, we conclude with
\cite[Lem.\ 2]{bj_Wagemann99} 
that $\iota^* \colon 
H^{2}_{\mathrm{ct}}(\Gamma(\ad(P))_{\mathbb{C}},W)
\rightarrow
H^2_{\mathrm{ct}}((\mathbb{C}[X]\otimes_{\mathbb{R}}\fg)^{\Delta},W)
$ is an isomorphism, 
where 
$W$ is a complex Fr\'{e}chet space considered as a trivial module,
and continuity is in the $C^{\infty}$-topology on \emph{both} sides.
%
Restricted to the equivariant map algebra, our canonical cocycle
takes values in the space $\ol{\Omega}{}_{\mathrm{alg}}^{1}(\mathbb{C}[X])$ of
K\"ahler differentials
modulo closed ones, and can be written
\begin{equation} \label{bj_alg_cocycle}
\omega_{\mathrm{alg}} : 
\wedge^{2} (\mathbb{C}[X]\otimes_{\mathbb{C}} \fg_{\mathbb{C}})^{\Delta} 
\to 
(\ol{\Omega}{}_{\mathrm{alg}}^{1}(\mathbb{C}[X])\otimes_{\mathbb{C}} 
V(\fg_{\mathbb{C}}))^{\Delta}  
\quad \colon \quad 
\xi \wedge \eta \mapsto [\kappa(\xi , d \eta)]\,.
\end{equation}
It is universal in the sense that
every \emph{continuous} $\mathbb{C}$-valued cocycle $\tau$ on the equivariant 
map algebra can be written up to coboundary as
$\tau = \phi \circ \omega_{\mathrm{alg}}$ for some continuous $\mathbb{C}$-linear  
functional $\phi$ on
$(\ol{\Omega}{}_{\mathrm{alg}}^{1}(\mathbb{C}[X])\otimes_{\mathbb{C}} 
V(\fg_{\mathbb{C}}))^{\Delta}$.

In the case of twisted multiloop algebras, a cocycle is continuous if it is of polynomial 
growth in the $\mathbb{Z}^{n}$-grading of $\mathbb{C}[T^{n}]$. 
If $\fg_{\mathbb{C}}$ is simple, 
then $\kappa$ is just the $\Aut(\fg_{\mathbb{C}})$-invariant
Killing form, so that
$V(\fg_{\mathbb{C}}) \simeq \mathbb{C}$ is a trivial $\Delta$-representation.
The universal cocycle thus takes values in 
$\ol{\Omega}{}^{1}_{\mathrm{alg}}(\mathbb{C}[T^{n}])^{\Delta}$,
in agreement with
the purely algebraic result \cite{bj_JieSun}. 
It might not be overly optimistic to hope 
for universality of (\ref{bj_alg_cocycle}) 
for equivariant map algebras with semisimple $\fg_{\mathbb{C}}$
in a more general context.


\end{talk}

\end{document}